\magnification=\magstep1
\input amstex
\documentstyle{amsppt}
\pagewidth{16 true cm}
\pageheight{23 true cm}
\topmatter
\TagsOnRight
\NoBlackBoxes
\rightheadtext { }
\title  Local  zeta functions and Newton polyhedra
\endtitle
\rightheadtext{ }
\leftheadtext{}
\author
 W. A. Zuniga - Galindo*
\endauthor
\address
\noindent Department of Mathematics and Computer Science, Barry University, 11300 N.E.
Second Avenue, 
Miami Shores, Florida 33161, USA
\endaddress
\email \noindent wzuniga\@mail.barry.edu
\endemail
\thanks
* This work was supported by  COLCIENCIAS, project  1241-05-111-97
\endthanks
\abstract
To a polynomial $f$ over a non-archimedean local field $K$ and  a 
cha\-rac\-ter $\chi$ of the group  of units of the valuation ring of $K$ one associates 
Igusa's local zeta function $Z(s,f,\chi)$. 
In this paper,  we study the local zeta function  $Z(s,f,\chi)$ associated to a non-degenerate
 polynomial $f$, by using an approach based on the p-adic stationary phase formula and N\'eron p-desingularization. 
 We give a small set of candidates for 
the poles of $Z(s,f,\chi)$ in terms of the Newton polyhedron $ \Gamma(f)$ of $f$. We also show that for almost all $\chi$, the local zeta  function $Z(s,f,\chi)$ is a polynomial in $q^{-s}$ whose degree is bounded by a constant independent of $\chi$. Our second result is a description of the largest pole of $Z(s,f, \chi_{\text{triv}})$ in terms of $ \Gamma(f)$ when the distance between $\Gamma(f)$ and the origin is at most one.

\endabstract
\endtopmatter

\document
\openup 3pt

\head
 {\bf 1. Introduction}
\endhead

 Let $K$ be a  non-archimedean  local field, and let $\Cal{O}_{K}$ be  the   ring of integers 
of $K$ and $\Cal{P}_K$ its maximal ideal. Let $\pi$ be  a  fixed uniformizing parameter of $K$,
 and let the residue field of $K$ be
   $\Bbb{F}_q$ the field with $q = p^{r}$ elements.  For  $x \in K$ we denote by  $v$  the
 valuation of $K$ such that $v(\pi) =1$,   $ |x|_{K} = q^{-v(x)}$ its absolute value and
 $ac(x) = x\pi ^{-v(x)}$ its angular component. 
Let  $f(x) \in  \Cal{O}_{K}[x]$, $x =(x_1,..,x_n)$ be a non-constant polynomial, and $ \chi : \Cal{O}_{K}^{\times} \longrightarrow \Bbb{C}^{\times}$ a 
character of $\Cal{O}_{K}^{\times}$, the group of units of $\Cal{O}_{K}$. We formally put
 $\chi(0) =0$. To these data one associates Igusa's local zeta function,

$$ Z(s,f,\chi)= \int_{\Cal{O}^{n}_{K}}\chi(ac f(x))|f(x)|_K ^{s} \mid dx \mid , \,\,\,\,\,\,\,\,  s \in \Bbb{C}, \,\,$$

\noindent for $ Re(s) >0$, where  $\mid dx \mid$ denotes the  Haar measure on $K^{n}$, normalized such that  
$\Cal{O}^{n}_{K}$ has measure 1. 
 In  the case
of  $K$ having  characteristic zero, Igusa [I2] and Denef [D1] proved that  $Z(s,f,\chi)$ is a rational  function of $q^{-s}$.

A basic problem is to determine the poles of the 
meromorphic  continuation  of 
$Z(s,f,\chi)$ 
into $Re(s) < 0$. 
 The general strategy  is to take  a   reso\-lution $h: X \longrightarrow K^{n}$
of $f  $ and study the  resolution data $\{ (N_i ,n_i)\}$ in which $N_i$ is the multiplicity
 of $f\circ h $ along a 
exceptional divisor $D_i$, and  $n_i$ is the multiplicity of $h^{\star}(dx)$ along $D_i$.  The set of 
ratios
$\{\frac{-n_i}{N_i}\} \bigcup \{-1\}$ contains the real parts of the poles of  $Z(s,f,\chi)$ as observed in 
[I2]. However, many examples show that most of these  ratios do not correspond to poles. The problem of the determination of  the actual poles  of  $Z(s,f,\chi)$ for arbitrary $n$ is still an open problem. The case $n=2$ was solved  for irreducible $f$ and $\chi= \chi_{\text{triv}}$ for all primes $p$ by Meuser [Me]. The generalization to reducible $f$ and  $\chi \ne \chi_{\text{triv}}$  but for almost all primes  $p$ was solved by Veys in [Ve].

In  case of non-degenerate polynomials with respect to its Newton polyhedron and  $K= \Bbb{R}$, Varchenko [Va] gave a procedure to compute a set of candidates  for  the poles of the complex power of $f$, by using  toroidal reso\-lution of singularities (see also [D-S-1], [D-S-2] for further generalizations). 

The p-adic case is entirely similar to the real case. In this case,  
Lichtin and Meuser [L-M] proved in the case $n=2$
that not all candidates provided by the numerical data of a toric resolution of $f$ are actually poles of $Z(s,f,\chi)$.
In [D3], Denef  gave  a procedure based on monomial changes of variables to determine  a small set of candidates  for the poles of
$ Z(s,f,\chi_{\text{triv}})$ in terms of the Newton polyhedron of $f$.

In this paper,  we study the local zeta function  $Z(s,f,\chi)$ associated to a globally non-degenerate
 polynomial $f$ (see  definition 1.1), by using an approach based on the p-adic stationary phase formula and N\'eron p-desingularization. We  show the stationary phase formula gives a small set of candidates for the poles  of $Z(s,f,\chi)$ in terms of the Newton polyhedron $ \Gamma(f)$ of $f$ (cf. theorem A). When $\chi= \chi_{\text{triv}}$ and char$(K) =0$ this set of poles agree with that obtained in [D3]. We also show that for almost all $\chi$, the zeta function $Z(s,f,\chi)$ is a polynomial in $q^{-s}$ whose degree is bounded by a constant independent of $\chi$. Our second   result  shows that the stationary  phase formula can be used to describe the largest pole of $Z(s,f, \chi_{\text{triv}})$  in terms  of $ \Gamma(f)$, when the distance  between $\Gamma(f)$ and the origin is at most one  (cf. theorem B). This result was previously known for char$(K)=0$.  This result allows one to generalize estimates for exponential sums that were obtained in [D-Sp] to the case char $(K) \ne 0$ (cf. corollary 6.1).

\smallskip
We set $\Bbb{R}_{+}= \{ x \in \Bbb{R} \,\, \mid \,\, x \geqq 0\}$. 
Let $f(x) = \sum_l a_l x^{l} \in K[x]$, $x= (x_1,x_2,...,x_n)$ be a polynomial in $n$ variables satisfying $f(0)=0$. 
 The set  $supp(f) = \{ l \in \Bbb{N}^{n} \,\ \mid \, a_l \ne 0\} $ is called the {\it support} 
of $f$. {\it The Newton polyhedron} $\Gamma(f)$ of $f$ is defined as the convex hull in 
$\Bbb{R}_{+}^{n}$ of the set

$$ \bigcup_{l \in supp(f)}\left( l +\Bbb{R}_{+}^{n} \right).$$
 We denote by $< , >$ the usual inner product of $\Bbb{R}^{n}$, and identify $\Bbb{R}^{n}$ with its dual  by means of it. We set 
$$<a_{\gamma},x> =m(a_{\gamma}),$$
 for the equation of the supporting hyperplane of a facet $\gamma$ (i.e. a face of codimension 1 of $\Gamma(f)$) with perpendicular vector
$a_{\gamma} =(a_{1},a_{2},...,a_{n}) \in \Bbb{N}^{n}\smallsetminus \{0\}$,
and   $\mid a_{\gamma}\mid:= \sum_i a_{i}$.

\definition {Definition 1.1}A polynomial  $f(x)=\sum_{i }a_i x^{i} \in K[x]$ is called  
 {\it globally non-degenerate  with respect to its Newton polyhedron } $\Gamma(f)$, if it satisfies the following two properties:

\noindent (GND1) the origin of $K^{n}$ is a singular point of $f(x)$;

\noindent (GND2) for  every  face $\gamma \subset \Gamma(f)$ (including $\Gamma(f)$ itself), the polynomial 
$$f_{\gamma}(x):= \sum_{i \in \gamma}a_i x^{i}$$ 
\noindent has the property that there is no $ x \in (K \smallsetminus\{0\})^{n}$ such that

$$f_{\gamma}(x) =\frac{\partial f_{\gamma}}{\partial x_1}(x)= ....= \frac{\partial f_{\gamma}}{\partial x_n}(x)
= 0.$$ 

\enddefinition

Our first   result  is  the following.  

\proclaim{ Theorem  A}
 Let  $K$ be a non-archimedean local field, and let $f(x) \in \Cal{O}_{K}[x]$ be  a polynomial  globally non-degenerate  with respect to its Newton polyhedron $\Gamma(f)$. 
Then  the Igusa local zeta function $Z(s,f,\chi)$ is a rational function 
of $q^{-s}$ satisfying:

 (i) if $s$ is a pole   of $Z(s,f,\chi)$, then 
$$ s= -\frac{\mid a_{\gamma} \mid}{m(a_{\gamma})} + \frac{2\pi i}{log \, q} \frac{k}{m(a_{\gamma})},
 \,\, k \in \Bbb{Z}$$
 for some facet  $\gamma$ of $\Gamma(f)$ with perpendicular $a_{\gamma}$, and $m(a_{\gamma}) \ne 0$,  or 
$$ s= -1 + \frac{2\pi i}{log \, q} k, \,\, k \in \Bbb{Z};$$

(ii) if $\chi \ne \chi_{\text{triv}}$ and  the order of $\chi$ does not divide any $m(a_{\gamma}) \ne 0$, where $\gamma$ is a facet  of $\Gamma(f)$, then  $Z(s,f,\chi)$ is a polynomial in  $q^{-s}$, and its degree is bounded by a constant independent of $\chi$.
 
\endproclaim

For a polynomial $f(x) \in K[x]$ globally non-degenerate  with respect to its Newton 
polyhedron $\Gamma(f)$, we set
$$
\beta(f) := \max_{\tau _j}\{-\frac{\mid a_j\mid}{m(a_j)}\},
$$
where $\tau _j$ runs through  all facets of  $\Gamma(f)$ satisfying $m(a_j) \ne 0$.  The point
 $$T_0 = (-\beta(f)^{-1},...,-\beta(f)^{-1})\in \Bbb{Q}^{n}$$ 
\noindent is the intersection point
 of the boundary of the Newton polyhedron $\Gamma(f)$ with the diagonal
 $\Delta = \{ (t,..,t) \; \mid \, t \in \Bbb{R} \}$ in $\Bbb{R}^{n}$. Let $\tau_0$ be 
the face of smallest dimension  of $\Gamma(f)$  containing  $T_0$, and   $\rho$ its codimension. 

If $g(x) \in \Cal{O}_K[x]$, $x=(x_1,..,x_n)$, we denote by  $\overline{g(x)}$ its reduction modulo 
$\Cal{P}_K$.

The second  result of this paper describes  the largest pole of $Z(s,f,\chi_{\text{triv}})$,
when  $\beta(f) \geq -1$.
\proclaim{Theorem B}
Let  $K$ be a non-archimedean local field, and let  $f(x) \in \Cal{O}_K [x]$ be  a globally  non-degenerate polynomial  with respect to 
its Newton polyhedron $\Gamma(f)$. If $\beta(f) > -1$, then $\beta(f)$ is a pole of $Z(s,f,\chi_{\text{triv}})$  of multiplicity $\rho$. If $\beta(f) =-1$, then $\beta(f)$ is a pole of $Z(s,f,\chi_{\text{triv}})$  of multiplicity  less than or equal to $\rho +1$. Moreover, 
if  every face $\gamma \supseteqq \tau_0$ satisfies  $\text{Card}\left(\{z\in \Bbb{F}^{\times \, n}_q  \mid  \bar{f_{\gamma}}(z)=0 \}\right)>0$, then the multiplicity of $\beta(f)$ is exactly $\rho +1$.

\endproclaim

The largest pole of  $Z(s,f,\chi_{\text{triv}})$ when $f$ is non-degenerate with respect to its Newton polyhedron $\Gamma(f)$ and $\beta(f)>-1$  follows  from observations made  by Varchenko in [Va] and was originally noted in the $p-$adic case in [L-M] (although it is misstated there as $\beta(f)\ne-1$). The  case $\beta(f)=-1$ is treated in [D-H]. The case  of $\beta(f)<-1$ is more difficult and is established in [D-H] with some additional conditions on $\tau_0$ by using a difficult result on exponential sums. Thus our Theorem B gives a different proof of the cases  where $\beta(f)\geqq -1$.

The organization of this paper is as follows.  In section 2, we review  Igusa's stationary phase formula. The results of this section generalize our previous results in [Z-G]. Section 3 contains
some basic results about Newton polyhedra. In section 4, we prove  theorem A. In section 5, we  prove theorem B. Section 6 contains some consequences of the main theorems. More precisely, we give  estimates for exponential sums involving  globally non-degenerate  polynomials (cf. corollary 6.1). In section 7, we compute explicitly the local zeta functions of some polynomials in two variables and discuss the relation between the largest pole of $Z(s,f, \chi_{\text{triv}})$ and $\beta(f)$.

\subheading {Acknowledgments}
I  wish to thank to  Jan Denef, Kathleen Hoornaert, and the referee for   their suggestions  which led to an improvement of this work.

\head
{\bf 2. Igusa's stationary phase  formula}
\endhead

In [I3] Igusa introduced the stationary phase formula for $\pi-$adic integrals  and 
suggested that a closer examination of this formula might lead to a new  proof of the rationality of 
$Z(s,f,\chi)$ in any characteristic. Following  this  suggestion  the author proved the rationality of the local zeta function $Z(s,f,\chi_{\text{triv}})$ attached to  a semiquasihomogeneous polynomial $f$ over an arbitrary non-archimedean local field [Z-G]. 

Let $L$  be a ring and $f(x) \in L[x]$, we denote by $V_f ( L)$ the corresponding $L-$hyper\-surface and by $Sing_f(L )$ the $L$-singular locus. 

We denote by $\overline{x}$ the image of an element of $\Cal{O}_K$ under the canonical 
homomorphism $ \Cal{O}_K \longrightarrow  \Cal{O}_K / \pi 
\Cal{O}_K \cong \Bbb{F}_q $, i.e. the reduction modulo  $\pi$. Given $f(x) 
 \in \Cal{O}_K[x] $ such that not all its coefficients are in
  $\pi \Cal{O}_K$,  we denote by $\overline{f(x)}$ the polynomial
   obtained by reducing modulo $\pi$  the coefficients of $f(x)$.

We fix a lifting $R$ of $\Bbb{F}_q$ in $\Cal{O}_K$. By definition, the   set $R$ is mapped bijectively onto
$\Bbb{F}_q$ by the canonical homomorphism $ \Cal{O}_K \longrightarrow  \Cal{O}_K / \pi 
\Cal{O}_K $. 

\noindent Let $f(x) \in \Cal{O}_K[x]$ be  a polynomial in $n$ variables,
  $ P_1= (y_{1},...,y_{n})\in
  \Cal{O}_K^{n}$,  and  $m_{P_1}=(m_1,...,m_n) \in \Bbb{N}^{n}$. We call a  $K^{n}-$isomorphism
$\Phi_{m_{P_1}}(x)$ a {\it dilatation}, if it has the form $\Phi_{m_{P_1}}(x)= (z_1,..,z_n)$, 
$ z_i = y_i+ \pi ^{m_{i}}x_i$, for each $i =1,2,..,n$. The  {\it  dilatation }of $f(x)$ at $P_1$ induced by $\Phi_{m_{P_1}}(x)$  is defined as 

$$f_{P_1}(x):= \pi^{-e_{P_1}}f(\Phi_{m_{P_1}}(x)),\tag{2.1}$$

\noindent where   $e_{P_1}$ is the minimum order   of $\pi$ in  the coefficients
of $f(\Phi_{m_{P_1}}(x))$.  We call the $K-$hypersurface $V_{f_{P_1}}(K)$ {\it the dilatation of }
$V_f(K)$ at $P_1$ induced by $\Phi_{m_{P_1}}(x)$; the number $e_{P_1}$  
{\it  the arithmetic multiplicity  of $f(x)$ at $P_1$} by $\Phi_{m_{P_1}}(x)$,
  and the set  $ S(f_{P_1})$, the lifting of $Sing_{\overline{f}_{P_1}}(\Bbb{F}_q)$, {\it the first generation of descendants
  of $P_1$}. 

 Given a sequence of dilatations $( \Phi_{m_{P_k}}(x)) _{k \in \Bbb{N}}$, we
 define inductively $e_{P_1,...,P_k}$ and $f_{P_1,...,P_k}(x)$, 
$ S(f_{P_1,...,P_k})$  as follows:
$$f_{P_1,...,P_k}(x):= \cases
f(x), \,\,\,\,\ \text{if} \,\,\ k=0,\\
\pi^{-e_{P_1,...,P_k}}f_{P_1,...,P_{k-1}}( \Phi _{m_{P_k}} (x)), \,\,\,\,\text{if} \,\,\,\, k\geqq 1,\\
\endcases
\tag{2.2}$$
\noindent where $P_k \in S(f_{P_1,...,P_{k-1}})$, and $e_{P_1,...,P_k}$  is the minimum order of $\pi$ in the coefficients of $f_{P_1,...,P_{k-1}}( \Phi _{m_{P_k}} (x))$. For  $k \geq 1$, the  set $S(f_{P_1,...,P_k}):=\bigcup_{P_k} S(f_{P_1,...,P_{k-1},P_k})$ is called {\it the $k^{th}-$generation of descendants of $P_1$}. By 
definition  the $0^{th}-$generation of descendants of $P_1$ is $\{P_1\}$.

\noindent Now, we review Igusa's stationary phase formula, from the point of view of the dilatations. For that, we  fix the $m_{P_k}$'s equal to $(1,..,1) \in \Bbb{N}^{n}$ in (2.1).

 Let $\overline{D}$ be a subset of $\Bbb{F}_q^{n}$ and 
  $D$ its preimage under the canonical homomorphism
 $\Cal{O}_K \longrightarrow  \Cal{O}_K / \pi 
\Cal{O}_K \cong \Bbb{F}_q $. Let  $S(f,D)$ denote   the subset of $R^{n}$
 (the set of representatives of $\Bbb{F}_q^{n}$ in $\Cal{O}_K^{n}$) mapped bijectively to the set $Sing_{\bar{f}}(\Bbb{F}_q) \bigcap \overline{D}$. We use the simplified notation $S(f)$ in the case of $D= \Cal{O}_K^{n}$. Also we define:
 
$$ \nu (\bar{f},D, \chi):= \cases
 q^{-n}\text{Card} \{\overline{P} \in \overline{D} |\,\,  
\overline{P} \notin V_{\bar{f}}(\Bbb{F}_q)\}, \,\, \text{if }\,\, \chi= \chi _{triv},\\ 
q^{-n c_{\chi}} \sum_{\{ P \in D \, \mid \, \overline{P} \notin V_{\bar{f}}(\Bbb{F}_q)\} \text{mod} \,\, \Cal{P}_K^{c_{\chi}} }\chi (ac(f(P))), 
 \,\, \text{if }\,\,\chi \ne\chi _{triv},\\
\endcases
$$
\smallskip

where    $c_{\chi}$ is  the conductor of $\chi$, and
$$\sigma (\bar{f},D,\chi) := \cases 
q^{-n}\text{Card} \{\overline{P} \in \overline{D} |\,\, 
\overline{P} \,\,\ \text{is a smooth point of} \,\,\,
 V_{\bar{f}}(\Bbb{F}_q) \}, \,\,  \text{if } \,\, \chi= \chi _{triv},\\ 
0, \,\,\, \text{if }\,\,\chi \ne\chi _{triv}.\\
\endcases
$$

 If $ D= \Cal{O}_K^{n}$, we use the simplified notation $\nu (\bar{f},\chi)$, $\sigma (\bar{f},\chi)$.
 We denote by $ Z(D,s,f,\chi)$ the integral $\int_{D}\chi (ac(f(x)))|f(x)|_{K} ^{s}\mid  dx \mid$. 
With all this, we are able to establish  Igusa's stationary phase formula for $\pi-$adic integrals ([I3], p. 177):

\subheading{ Igusa's Stationary Phase Formula}

$$ Z(D,s,f,\chi)= \nu(\bar{f},D, \chi) +\sigma(\bar{f},D,\chi) \frac{(1-q^{-1})q^{-s}}{(1-
q^{-1-s})} +$$
$$ \sum_{P \in S(f,D)} q^{-n -
e_{P}s}\int_{\Cal{O}_K^{n}} \chi (ac (f_P(x))|f_{P}(x)|_K ^{s}\mid dx \mid, \tag{2.3}$$

\noindent where $Re(s) >0$. The proof given by Igusa in [I3], for the case $\chi = \chi_{triv}$,
 generalizes literally to  arbitrary characters.

In [Z-G] the author  introduced the following index  of singularity
 at a point $P \in \Cal{O}_K ^{n}$,  satisfying $P \notin Sing_{f}(\Cal{O}_K)$.

\definition{Definition 2.1} Let $f(x)\in \Cal{O}_K[x]$ be
 a polynomial and $ P=(a_1,..,a_n) \in  \Cal{O}_K ^{n}$, such that $ P \notin  Sing_{f}(\Cal{O}_K)$. 
We define

$$L(f,P) := \text{Inf}
 \left( v(f(P)), v( \frac{\partial f}{\partial x_1}(P)),..,
 v( \frac{\partial f}{\partial x_n}(P))\right).$$
\enddefinition

\noindent It follows from the definition that   $L(f,P) =0$ if and only if  the polynomial
$$  \overline{f(x)} = \alpha_0 +
\sum_{j} \alpha_j ( x_j - \overline{a_j}) + ( \text{ degree }
 \geq 2) \in \Bbb{F}_q[x],$$
\noindent  satisfies  $\alpha_j \in  \Bbb{F}_q^{*}$ for some $j=0,1,2,..,n$.
 
The index $L(f,P)$ appears naturally associated to  Igusa's stationary phase, as 
it was already noted in [Z-G]. In addition, this  index   
plays an important role in the  construction  of the N\'eron $\pi-$adic  desingularization of the special fiber 
of smooth  schemes over $Spec(\Cal{O}_K)$ (see [A], [N]).

If $A \subseteq \Cal{O}_K ^{n}$, we denote by $A^{c}$ the complement of $A$ with respect to  $ \Cal{O}_K^{n}$.

  \proclaim{Proposition  2.2}
   Let  $D \subseteq \Cal{O}_K ^{n}$ be an open and compact subset, and  let $f(x)\in \Cal{O}_K[x]$ be a  polynomial  such that
    $Sing_{f}(K) \bigcap D = \emptyset$.
Then there exists a   constant $C(f,D)\in \Bbb{N}$, depending only on $f$ and $D$, such that
$$ L(f,P) \leqq C(f,D), \,\,\, \text{for  all} \,\,\ P \in  D.\tag{2.4}$$
\endproclaim
\demo{Proof}
By contradiction, we suppose that $L(f,P)$ is not bounded on $D$. Thus there exists a sequence
$(Q_i)_{i \in \Bbb{N}}$ of points of $D$ satisfying  lim $L(f,Q_i) \,\, \longrightarrow \infty$, 
when $ i \longrightarrow \infty$. This sequence has a limit point $Q_{*} \in D$. 
Since $ Sing_{f}(K) $ is a closed set, we  have that $Q_{*} \in  Sing_{f}(K) \bigcap D = \emptyset$,  contradiction.
\quad\qed
\enddemo  
From now on, we shall suppose  that  $C(f,D)$ is   minimal for  condition  (2.4).

We recall  that a subset $A$ of  $K^{n}$ is open and compact if and only if  there 
is $m \geqq 0$  such that  $A$ is the finite union of classes modulo $\pi ^{m}$. In particular the preimage of any subset of $\Bbb{F}_q^{n}  $ under the canonical homomorphism  $\Cal{O}_K \longrightarrow  \Cal{O}_K / \pi \Cal{O}_K $ is an open and compact subset.

The  following lemma is a  generalization  of  proposition 2.3  of [Z-G].

\proclaim{Lemma 2.3}
 Let  $D \subseteq \Cal{O}_K^{n}$ be
the preimage under the canonical homomorphism 
$\Cal{O}_K \longrightarrow  \Cal{O}_K / \pi 
\Cal{O}_K $ of a subset $\overline{D} \subseteq \Bbb{F}_q^{n} $, and let   
 $f(x) \in \Cal{O}_K[x] $ be a
 polynomial such that
    $Sing_{f}(\Cal{O}_K) \bigcap D = \emptyset$, then

\noindent (i) $ L(f_{P_1,..,P_k},0) \leqq L(f,P_1 + \pi P_2+..+\pi ^{k-1}P_k) -k$, for every 
$P_k$, $k \geq 1$,  satisfying: (H1) $P_k$ is in the $(k-1)^{th}$-generation  of descendants of $P_1$; (H2) $P_k$ has
at least one descendant in the $k^{th}$-generation of descendants of $P_1$.

\noindent (ii) For any $P=P_1 \in S(f,D)$, if $k \geq C(f,D) +1$ then   $S(f_{P_1,P_2,..,P_k})= \emptyset$.
  
\endproclaim
\demo{Proof}
First, we observe that 
$$f(P_1+\pi P_2+..+\pi ^{k-1}P_k +\pi^{k}x)= 
\pi ^{E(P_1,..,P_k)}f_{P_1,..,P_k}(x),\tag{2.5}$$
where $E(P_1,..,P_k) = e_{P_1}+e_{P_1,P_2} +e_{P_1,..,P_k}$. The result follows from
(2.5), if 
 $$ e_{P_1,..,P_l} \geq  2, \,\, \text{ for} \,\, l=1,2,..,k.$$
 This last fact  follows from the following reasoning. 

 By applying the Taylor formula to $f_{P_1,..,P_{l-1}}(P_l +\pi x)$, we obtain
 $$ f_{P_1,..,P_{l-1}}(P_l +\pi x) = f_{P_1,..,P_{l-1}}(P_l) +\pi \sum_j
 \frac{\partial f_{P_1,..,P_{l-1}}}{ \partial x_j }(P_l) x_j +\pi ^{2} (\text{degree} \geq 2).\tag{2.6}$$
 From hypothesis (H1) follows that  $v( f_{P_1,..,P_{l-1}}(P_l)) \geq 1$ and $v(\frac{\partial f_{P_1,..,P_{l-1}}}{ \partial x_j }(P_l)) \geq 1$, and from hypothesis (H1) and (H2) that
 $$ v (f_{P_1,..,P_{l-1}}(P_l)) \geq 2;$$
 therefore  (2.6) implies that
$ e_{P_1,..,P_l}\geq 2$, $l=1,2,..,k$.

(ii) The second part of the lemma follows immediately from  (i).
\quad\qed
\enddemo

We observe  that if $P_l \in S(f_{P_1,..,P_{l-1}})$ does not have   descendants  in the $l^{th}-$generation (i.e. $S(f_{P_1,..,P_{l-1}, P_l})=\emptyset$), then the polynomial 
$$ f_{P_1,..,P_{l-1},P_l}(P_{l+1}+ \pi x) =   f_{P_1,..,P_{l}}(P_{l+1}) +\pi\sum_j
 \frac{\partial f_{P_1,..,P_{l}}}{ \partial x_j }(P_{l+1}) x_j +\pi^{2}(\text{degree} \geq 2)$$
satisfies 
$\overline{f_{P_1,..,P_{l}}(P_{l+1}) }\ne 0$, or $\overline{ \frac{\partial f_{P_1,..,P_{l}}}{ \partial x_{j_0} }(P_{l+1})} \ne 0$,
   for some $j_0$.  Thus  for any
$P_{l+1}$ satisfying 
$\overline{f_{P_1,..,P_{l}}(P_{l+1})} = 0$,  it holds  that $\overline{f_{P_1,..,P_{l+1}}(x)} $ is a polynomial of degree at most one.

\proclaim{Lemma 2.4}
 Let  $D \subseteq \Cal{O}_K^{n}$ be
the preimage under the canonical homomorphism $\Cal{O}_K \longrightarrow  \Cal{O}_K / \pi 
\Cal{O}_K $ of a subset $\overline{D} \subseteq \Bbb{F}_q^{n} $. Let   
 $f(x) \in \Cal{O}_K[x] $ be a
 polynomial such that
    $Sing_{f}(K) \bigcap D = \emptyset$, then
$$
\int_{D}\chi(ac f(x)) |f(x)|_K ^{s}|dx| =
\cases
\frac{T(q^{-s})}
{1-q^{-1}q^{-s}},
 \,\,\, 
\chi = 
\chi _{\text{triv}},
\\
L(q^{-s}), \,\,  \chi \ne \chi _{\text{triv
}},\\
\endcases
$$
\noindent where $T$ and $L$ are polynomials in $q^{-s}$ with rational coefficients. Furthermore, in the case $ \chi \ne \chi_{\text{triv}}$, the degree of the polynomial $L(q^{-s})$ is bounded by a constant depending only on $f$ and $D$.
\endproclaim

\demo{Proof}   
We define inductively $I_k$ as follows:

$$ I_1 := S(f,D),$$
$$ I_k :=\{(P_1,P_2,..,P_k) \,\, \mid \, (P_1,P_2,..,P_{k-1}) \in I_{k-1}, \, \text{and } \,
P_{k}\in S(f_{P_1,P_2,..,P_{k-1}})	\},\,\ k \geq 2.$$
We set $E(P_1,...P_k):= e_{P_1}+e_{P_1,P_2}+...+e_{P_1,P_2,..,P_k}$.

If $m= C(f,D)+1$,  then   $I_{m+1}=\emptyset$, because lemma 2.3 (ii) implies that
$S(f_{P_1,P_2,..,P_m})=\emptyset$, for every $(P_1,P_2,..,P_m) \in I_m$. By  applying the stationary phase formula  $m+1-$times, we obtain 
\bigskip

$$
\align
Z(D,s,f,\chi) = \nu (\bar{f},D, \chi)+
 \sigma (\bar{f},D, \chi)\frac{(1-q^{-1})q^{-s}}{(1-q^{-1-s})}+\\
\sum_{k= 1}^{m}q^{-kn}\left( 
\sum_{(P_1,...P_k) \in I_k}\nu(\overline{f}_{P_1,..,P_k},\chi)q^{-E(P_1,..,P_k)s}\right)+\\
 \frac{(1-q^{-1})q^{-s}}{(1-q^{-1-s})} \sum_{k = 1}^{m}q^{-kn}\left( 
\sum_{(P_1,...P_k) \in I_k}\sigma(\overline{f}_{P_1,..,P_k},\chi)q^{-E(P_1,..,P_k)s}\right).\\
\tag{2.7}
\endalign
$$

\noindent In the case  $ \chi \ne \chi_{triv}$, all 
$ \sigma(\overline{f}_{P_1,...P_k},\chi)=0$, thus 
$Z(D,s,f,\chi)$ is a polynomial in $q^{-s}$ and its degree is bounded by the maximum of 
the $E(P_1,...P_{m})$, where $P_m$ runs through the descendants of the 
$C(f,D)+1-$generation  of $S(f,D)$. 
\quad\qed 
\enddemo

\proclaim{Corollary 2.5}
 Let  $D \subseteq \Cal{O}_K^{n}$ be
the preimage under the canonical homomorphism $\Cal{O}_K \longrightarrow  \Cal{O}_K / \pi 
\Cal{O}_K $ of a subset $\overline{D} \subseteq \Bbb{F}_q^{n} $. Let   
 $F(x)= f(x)+\pi ^{\beta}g(x) \in \Cal{O}_K[x] $ be a
 polynomial such that $\beta \geqq C(f,D)+1$, and
$$Sing_{F}(K) \bigcap D =Sing_{f}(K) \bigcap D = \emptyset. $$
\noindent Then $$Z(D,s,F,\chi) =Z(D,s,f,\chi).\tag{2.8}$$

\endproclaim
\demo{Proof}
The result follows immediately from  expansion 2.7 and  the fact that $C(f,D)= C(F,D)$.
\quad\qed
\enddemo

\head
{\bf 3. Newton polyhedra}
\endhead
In this section we review some well-known results about Newton polyhedra that we shall use in this paper  (see  e.g. [K-M-S], [D-3]).

 We set $\Bbb{R}_{+}= \{ x \in \Bbb{R} \,\, \mid \,\, x \geqq 0\}$. 
 Let $f(x) = \sum_l a_l x^{l} \in K[x]$, $x= (x_1,x_2,...,x_n)$ be a 
 polynomial in $n$ variables satisfying  $f(0)=0$.  The set  $supp(f) = \{ l \in \Bbb{N}^{n} \,\ \mid \, a_l \ne 0\} $ 
 is called the {\it support} of $f$. {\it The Newton polyhedron} $\Gamma(f)$ of $f$ is defined as the convex hull in $\Bbb{R}_{+}^{n}$ of the set

$$ \bigcup_{l \in supp(f)}\left( l +\Bbb{R}_{+}^{n} \right).$$
{
By a {\it proper face} $\gamma$ of $\Gamma(f)$, we mean the non-empty convex set
$\gamma$ obtained by intersecting $\Gamma(f)$ with an affine hyperplane $H$, 
 such that $\Gamma(f)$ is contained in  one of two half-spaces determined
 by $H$. The hyperplane $H$ is named {\it the supporting hyperplane} of $\gamma$. 
 A face of codimension one is named a {\it facet}.

\noindent We set $<\, ,>$ for the usual inner product in  $\Bbb{R}^{n}$, and identify the dual vector space with $\Bbb{R}^{n}$. For $a \in \Bbb{R}_{+}^{n} $, we define
 
$$m(a) :=\inf_{x \in \Gamma(f)}  \{ <a,x> \}.$$
\noindent {\it The first meet locus} of $a  \in \Bbb{R}_{+}^{n} \smallsetminus \{0\}$ is defined by
$$ F(a): = \{ x \in \Gamma(f)\,\ \mid \, <a,x> = m(a)  \}.$$
  
The first meet locus $F(a)$ of $a$ is a  proper face of $\Gamma(f)$. 

We define an equivalence relation on $ \Bbb{R}_{+}^{n}\smallsetminus\{0\} $ by 
$$ a \backsimeq a' \,\,\, \text{if and only  if}\,\, F(a) =F(a').$$

If $\gamma$ is a face  of $\Gamma(f)$, we define the cone associated to $\gamma$ as
$$\Delta_{\gamma} := \{ a  \in (\Bbb{R}_{+ })^{n} \smallsetminus \{0\}  \, \mid \, F(a)=\gamma  \}.$$

\noindent The following  two propositions  describe the geometry of the equivalences classes  of $\backsimeq$ (see e.g. [D-3]).

\proclaim{Proposition 3.2}
 Let $\gamma$ be a proper face of $\Gamma(f)$. Let $ w_1,w_2,..,w_e$ be the facets of $\Gamma(f)$ which contain $\gamma$. Let $a_1,a_2,..,a_e$ be vectors which are perpendicular to respectively $w_1,w_2,..,w_e$. Then 

$$\Delta_{\gamma}= \{ \sum_{i=1}^{e} \alpha_i a_i \,\, \mid \, \alpha_i \in \Bbb{R}, \,\ \alpha_i > 0\}.$$
\endproclaim

If $a_1,a_2,..,a_e \in \Bbb{R}^{n}$, we call $ \{ \sum_{i=1}^{e} \alpha_i a_i \,\, \mid \, \alpha_i \in \Bbb{R}, \,\ \alpha_i > 0\} $ {\it the cone strictly positive spanned} by the vectors $a_1,a_2,..,a_e$. Let $\Delta$ be a cone strictly positive spanned by the vectors $a_1,a_2,..,a_e$.  If $a_1,a_2,..,a_e$ are linearly  independent over $\Bbb{R}$, the cone $\Delta$ is called {\it  a simplicial cone}. In this last case, if $a_1,a_2,..,a_e \in \Bbb{Z}^{n}$, the cone $\Delta $ is called {\ a rational simplicial cone}.
If $\{a_1,a_2,..,a_e\}$ can be completed to be 
 a basis  of $\Bbb{Z}-$module $ \Bbb{Z}^{n}$, the cone $\Delta$ is named {\it  a simple cone}.

A vector $ a \in \Bbb{R}^{n}$  is called {\it primitive}  if the components of $a$ are positive integers whose greatest common divisor is one.

For every facet of $\Gamma(f)$  there is a unique primitive vector in $\Bbb{R}^{n}$ which is perpendicular to this facet. Let $ \Cal{D} $ be the set of all these vectors.

\proclaim{Proposition  3.3}
Let $\Delta$ be the cone strictly positively spanned by vectors $a_1,a_2$,..,$a_e \in \Bbb{R}_{+}^{n} \smallsetminus \{0\}$. Then there is a partition of $\Delta$ into cones $\Delta_i$, such that each $\Delta_i$ is strictly positively spanned by some vectors from $\{ a_1,a_2,..,a_e\}$ which are linearly independent over $\Bbb{R}$.
\endproclaim

The  two previous  propositions imply  the  existence of  a partition of $\Delta_{\gamma}$ into  rational simplicial cones.

\proclaim{Proposition 3.4}([K-M-S], p. 32-33) Let $\Delta$ be a rational simplicial cone. Then there exists a partition of $\Delta$ into simple cones.
\endproclaim
 
Summarizing, given a polynomial $f(x) \in K[x]$, $f(0)=0$,  with Newton polyhedron $\Gamma(f)$, there exists a finite  partition of $\Bbb{R}_{+}^{n}$  of the form:

$$\Bbb{R}_{+}^{n} = \{(0,..,0)\}\bigcup  \bigcup_i \Delta_i,$$
 
\noindent where each $\Delta_i$ is  a simplicial  cone contained in an equivalence class of $\backsimeq$. Furthermore, by proposition 3.4, it is possible to refine this partition in such a way that each $\Delta_i$ is  a simple  cone contained in an equivalence class of $\backsimeq$.

\head
{\bf 4. Local zeta functions of  globally non-degenerate polynomials}
\endhead
 In this section we  prove   theorem A. First, we give some preliminary results.

If $A \subseteq \Bbb{Z}_{+}^{n}$, we  set 
$$E_{A}:=\{(x_1,..,x_n) \in \Cal{O}_{K}^{n}\, \mid \,\, (v(x_1),..,v(x_n)) \in A\},$$
and  
$$
Z_{A}(s,f,\chi) :=\int_{E_{A}} {\chi(ac f(x))|f(x)|_K ^{s} \mid dx \mid}.
$$
 Also, if $B\subseteq \Cal{O}_{K}^{n}$,  we set 
$$
Z(B,s,f,\chi) :=\int_{B} {\chi(ac f(x))|f(x)|_K ^{s} \mid dx \mid}.
$$
 Thus $Z_{A}(s,f,\chi)=Z(E_A,s,f,\chi)$.

\proclaim{Proposition 4.1}
Let $f(x) \in \Cal{O}_K[x]$ be a globally non-degenerate polynomial with respect to its 
Newton polyhedron $\Gamma(f)$,  $\gamma \subseteq \Gamma(f)$   a proper face,  and  $\Delta_{\gamma}$ its associated  cone. If  $\Delta_{\gamma}$ is a simple cone  spanned by $a_1, a_2,..,a_e \in \Cal{D}$, and $f(x) = f_{\gamma}(x) +\pi ^{g_0}H(x)$, where  $g_0 \geq C (f_{\gamma}, \Cal{O}_K^{\times})+1$ (the  constant whose existence  was established in proposition  2.2), and  all monomials of $H(x)$ are not in $\gamma$, then 
$$
Z_{\Delta_{\gamma} }( s,f, \chi )= Z(\Cal{O}_K ^{\times n }, s,f_{\gamma}, \chi )
\frac{q ^{-\Sigma_{j=1}^{e}(\mid a_j\mid +m(a_j)s)}}
{ \prod_{j =1}^{e} (1- q ^{-\mid a_j\mid -m(a_j) s})}.
\tag{4.1}
$$
\endproclaim
\demo{Proof}
The hypothesis $\Delta_{\gamma} $ is a simple cone spanned by $a_j=(a_{1,j},a_{2,j},..,a_{n,j})$, $j=1,2,..,e$, implies that
 $$\Delta_{\gamma} \bigcap \Bbb{N}^{n} = \bigoplus_{j=1}^{e}a_j(\Bbb{N} \smallsetminus \{0\}). \tag{4.2}$$

From (4.2), we obtain the following expansion for $ Z_{\Delta}( s,f, \chi )$:

$$
Z_{\Delta_{\gamma}}( s,f, \chi )= \sum_ {y_1=1}^{\infty}..\sum_{y_e = 1}^{\infty} \int_{\omega_ {(y_1,.., y_e)}}
\chi(ac f(x))
\mid f(x)\mid _K ^{s}
\mid dx \mid, 
\tag{4.3}
$$
 where  
$$\omega_ {(y_1,..,y_e)}:= \{(x_1,..,x_n)  \in \Cal{O}_K ^{n} \,\ \mid \,\ x_i =
\pi^{\Sigma_j a_{i,j}y_j}
\mu_i, \,\, \mu_i \in \Cal{O}_K^{\times}, \,\, i=1,2,..,n\}.$$

In order to compute  the integral in (4.3), we introduce the  dilatation
 $$ \Phi_{(y_1,..,y_e)}(x)= (\Phi_1(x),..,\Phi_n(x)) : K ^{n} \longrightarrow K^{n},$$
 where
$$ \Phi_i ( x)=   \pi^{\Sigma_j a_{i,j}y_j} x_i, \,\, i=1,2,..,n. \tag{4.4}$$ 

 By using the dilatation  (4.4) as a change of variables in (4.3), it holds that
$$
\int_{\omega_ {(y_1,..,y_e )}}\chi(ac f(x))\mid f(x)\mid _K ^{s}\mid dx \mid= $$
$$
q ^{-\Sigma_{j=1}^{e}y_j(\mid a_j\mid +m(a_j)s)}
\left(\int_{\Cal{O}_K^{\times n}}\chi(ac (f_{(y_1,..,y_e )}(x)) 
\mid f_{(y_1,..,y_e )}(x)\mid _K ^{s}\mid dx \mid \right),
\tag{4.5}
$$
where $f_{(y_1,..,y_e)}(x)= f_{\gamma}(x) +\pi^{g(y_1,..,y_e)+g_0}H_{(y_1,..,y_e)}(x)$,  and $g(y_1,..,y_e) \geq 1$. 
The result follows from (4.5) by using corollary 2.5  and  expansion (4.3).\quad\qed 
\enddemo

\proclaim{Proposition 4.2}
Let $f(x) \in \Cal{O}_K[x]$ be a   globally non-degenerate polynomial  with respect to its Newton
polyhedron $\Gamma(f)$,  $\gamma \subseteq \Gamma(f)$   a proper face,  and  $\Delta_{\gamma}$ its associated  cone. If $\Delta_{\gamma}$ is a simple cone  spanned by $a_1, a_2,..,a_e \in \Cal{D}$, then 
$$
\align
Z_{\Delta_{\gamma}}( s,f, \chi )= 
\sum_{ y \,\, \text{finite}} A_y(q^{-s})Z(\Cal{O}_K ^{\times n }, s,f_y, \chi )+
\sum_{I \subseteq \{1,2,..,e\}}\frac{A_I(q ^{-s})Z(\Cal{O}_K ^{\times n }, s,f_I, \chi )}
{ \prod_{j \in I}(1- q ^{-\mid a_j\mid -m(a_j) s})},\\
\endalign
$$
where $y$  runs through a finite number of points in $ \Bbb{N}^{n}$,  $A_y(q^{-s})$, $A_I(q^{-s}) \in \Bbb{Q}[q^{-s}]$,   $f_y(x)$ and  $f_I(x)$ are polynomials in $\Cal{O}_K[x]$  satisfying $Sing_{f_y}(K) \bigcap (K \smallsetminus \{0\})^{n} = \emptyset$, for every $y \in \Bbb{N}$,  $Sing_{f_I}(K) \bigcap (K \smallsetminus \{0\})^{n} = \emptyset$, for every $I$, respectively. Furthermore, if $\gamma_{a_i}$ denotes the facet with perpendicular $a_i$, and $\gamma_I = \bigcap_{i \in I} \gamma_{a_i}$, then $f_I(x)= f_{\gamma_{I}}(x)$.
\endproclaim
\demo{Proof}
By induction on $l$, the number of generators of  the simple cone $\Delta_{\gamma}$.

{\bf Case} l=1

Let $m_0=   C(f_{\gamma}, \Cal{O}_K^{\times})+1$, and
$$S := \Delta_{\gamma} \bigcap \Bbb{N}^{n}= \{a_1 y \,\mid \, y\in \Bbb{N}, y\geq 1\}.$$
The set $S$ can be partitioned into the subsets $S_0$, $S_1$, defined as  follows:
$$S_0 :=\{a_1 y \,\mid \, y=1,2,..,m_0-1\}, \,\,\, 
 S_1:= \{a_1 y \,\mid \,y \in \Bbb{N}, y \geqq m_0\}.$$
Also we define
$$E_0:= \{(x_1,..,x_n)\in \Cal{O}_K^{n} \, \mid \,(v(x_1),..,v(x_n))\in S_0\},$$
$$E_1:= \{(x_1,..,x_n)\in \Cal{O}_K^{n} \, \mid \,(v(x_1),..,v(x_n))\in S_1\}.$$

Thus $ 
Z_{\Delta_{\gamma}}( s,f, \chi )=  Z(E_0,s,f, \chi) + Z(E_1,s,f,\chi) $, and by making a change of variables of type (4.4), we obtain

$$
\align
Z_{\Delta_{\gamma}}( s,f, \chi )= \sum_{y=1}^{m_0-1}q^{-y ( \mid a_1 \mid+m(a_1)s)} Z(\Cal{O}_K^{\times},s,f_y,\chi) +\\
q^{-m_0 ( \mid a_1 \mid + m(a_1)s)}Z_{\Delta_{\gamma}}(s,f_{a_1}(x)+ \pi ^{m_0 }H(x),\chi),\\
\tag{4.6}
\endalign$$

\noindent where $f_y(x)$ are obtained from $f(x)$ by a change of variables of type (4.4) followed  by a division  by a power of $\pi$, $f_{a_1}(x)$ is the restriction of $f(x)$ to the facet $\gamma_{a_1}$ with perpendicular $a_1$, and all monomials of $H(x)$ are not  in  $\gamma_{a_1}$. The result follows from (4.6),  by means of the following  equality (cf. proposition 4.1)
$$q^{-m_0 ( \mid a_1 \mid + m(a_1)s)}Z_{\Delta_{\gamma}}(s,f_{a_1}(x)+ \pi ^{m_0 }H(x),\chi)=
\frac{q^{-(m_0 +1)( \mid a_1 \mid + m(a_1)s)}}{1-q^{-( \mid a_1 \mid + m(a_1)s)}}Z(\Cal{O}_K^{\times},s,f_{a_1},\chi).$$

{\bf Induction hypothesis}
Suppose that the lemma is valid for every   polynomial $f(x)$ globally non-degenerate with respect its Newton polyhedron, and  for every  simple cone  spanned by  at most $e-1$ vectors of $\Cal{D}$.

{\bf Case $l >1$}

Let $f(x)$ be  globally non-degenerate  polynomial  and $\Delta_{\gamma}$  a simple cone spanned by $a_1, a_2,..,a_e$,  satisfying the conditions of proposition 4.2.

 We set  $m_0= C( f_{\gamma},   \Cal{O}_K^{\times})+1$, and

$$S := \Delta_{\gamma} \bigcap \Bbb{N}^{n}= \bigoplus _{j=1}^{e}a_j (\Bbb{N}\smallsetminus \{0\}),\tag{4.7}$$
 $a_j = (a_{1,j},..,a_{n,j})$, $j=1,2,..,e$. For each subset $I \subseteq \{1,2,..,e\}$, we put 
 $r_I \in \Bbb{N}^{e -\text{Card}(I)}$, $r_I =(r_{i_1},r_{i_2},..,r_{i_{e -\text{Card}(I)}})$, with $0 <  r_{i_l} \le m_0-1$, $l=1,2,..,e -\text{Card}(I)$.
The set $S$  admits the  following partition:
$$ S = \bigcup_{I,r_{I}}S_{I,r_I},\tag{4.8}$$
with
$$S_{I,r_I} = \{ \sum_{j \in I}a_j y_j + \sum_{j \notin I} a_j r_j  \mid \,  y_j \geq m_0, \,\, \text{if}
\,\,  j \in I,  \, \text{and} \,\,   y_j = r_{i_j} \,\, ,  \text{if} \,\, j \notin I \},
$$ where for each $I \subseteq \{1,2,..,e\}$, the corresponding  $r_I$'s run through all possible different integer vectors satisfying the above mentioned conditions.
We set 
$$E_{I,r_I}:=\{(x_1,..,x_n) \in \Cal{O}_{K}^{n}\, \mid \,\, (v(x_1),..,v(x_n)) \in S_{I,r_I}\}.$$

It follows from  partition (4.8) that 
$$Z_{\Delta_{\gamma}}(s,f,\chi) =\sum_{I,r_I}Z(E_{I, r_I},s,f,\chi). \tag{4.9}$$

By a change of variables of type 
$$ \Phi_i(x)= \pi^{(\sum_{j \in I}a_{i,j} y_j + \sum_{j \notin I}a_{i,j} r_j)}x_i,\,\,  i=1,..,n;$$ the integral $Z(E_{I, r_I},s,f,\chi)$ equals
$$q^{-m_0\Sigma_{j\in I}(\mid a_j\mid + m(a_j)s)-\Sigma_{j\notin I}r_j(\mid a_j\mid + m(a_j)s)}Z_{\Delta_I}(s,f_I,\chi),\tag{4.10}$$
where $\Delta_I$ is 
 a simple cone generated by $a_i$, $i \in I$, and $f_I(x)$ is obtained from  $f(\Phi_i(x))$ by  division by a power of $\pi$. From these observations and (4.9), we obtain

$$
Z_{\Delta_{\gamma}}(s,f,\chi)= \Sigma_{ I \subset \{1,2,..,e\}} A_{I}(q^{-s}) Z_{\Delta_I}(s,f_I,\chi)+
$$

$$
q^{-m_0 \Sigma_{j=1}^{e}(\mid a_j\mid +m(a_j)s)}Z_{\Delta_{\gamma}}(s,f_{\gamma}+ \pi^{g_0}H(x),\chi),
\tag{4.11}
$$
where $I$ runs through all proper subsets of $\{1,2,..,e\}$,  $A_I(q^{-s})=\sum_k q^{-a_k(I)-b_k(I)s}$, $a_k(I), b_k(I) \in \Bbb{N}$,   $ g_0 \geqq m_0$, and all monomials of $H(x)$ are not in  $\gamma$.  From (4.11) and proposition 4.1, we obtain
$$
Z_{\Delta_{\gamma}}(s,f,\chi)= \Sigma_{ I \subset \{1,2,..,e\}} A_{I}(q^{-s}) Z_{\Delta_I}(s,f_I,\chi)+
$$

$$
q^{-(1+m_0) \Sigma_{i=1}^{e}(\mid a_i\mid +m(a_i)s)}Z(\Cal{O}_K ^{\times n }, s,f_{\gamma}, \chi )
\frac{1}
{ \prod_{j =1}^{e} (1- q ^{-\mid a_j\mid -m(a_j) s})}.
\tag{4.12}
$$
 
The result follows from the induction hypothesis and (4.12).\quad\qed
\enddemo

  We observe that each  $A_I(q^{-s})$ in  proposition 4.1  is a finite sum of monomials of type $q^{-a_I-b_I s}$, with  $a_I$, $b_I >0$. We also note that
a facet with supporting hyperplane   $ x_{i_0} =0$ contributes to the denominator of $Z_{\Delta_{\gamma}}(s,f,\chi)$ with a constant factor  $\frac{1}{1-q^{-1}}$. 

\bigskip
The proof of proposition 4.2 can be  easily adapted to state the following more general result.

\proclaim{Corollary 4.3}
Let $f(x) \in \Cal{O}_K[x]$ be a   globally non-degenerate polynomial with respect to its Newton polyhedron $\Gamma(f)$,  $\gamma \subseteq \Gamma(f)$   a proper face, and    $\Delta_{\gamma}$ its associated   cone. Let    $\{a_1,a_2,..,a_f\} \subset \Cal{D}$ be a  set of generators of  $\Delta_{\gamma}$,  $\{a_1,a_2,..,a_e\} \subset \{a_1,a_2,..,a_f\}$ of  $e$ $\Bbb{R}-$linearly independent vectors, and $b \in \Delta_{\gamma} \bigcap \  (\Bbb{N}\smallsetminus \{0\})^{n}$. We set $\Delta:= b +  \bigoplus _{j=1}^{e}a_j \Bbb{N}$. Then 

$$
Z_{\Delta}(s,f,\chi)=
\sum_{ y } A_y(q^{-s})Z(\Cal{O}_K ^{\times n }, s,f_y, \chi )+
\sum_{I \subseteq \{1,2,..,e\}}\frac{A_I(q ^{-s})Z(\Cal{O}_K ^{\times n }, s,f_I, \chi )}
{ \prod_{j \in I}(1- q ^{-\mid a_j\mid -m(a_j) s})},
$$
where $y$  runs through a finite number of points in $ \Bbb{N}^{n}$,  $A_y(q^{-s})$, $A_I(q^{-s}) \in \Bbb{Q}[q^{-s}]$,
with  $A_I(q^{-s})=\sum_k q^{-a_k(I)-b_k(I)s}$, $a_k(I), b_k(I) \in \Bbb{N}$,  $f_y(x)$ and  $f_I(x)$ are polynomials in $\Cal{O}_K[x]$  satisfying $Sing_{f_y}(K) \bigcap (K \smallsetminus \{0\})^{n} = \emptyset$, for every $y$,  $Sing_{f_I}(K) \bigcap (K \smallsetminus \{0\})^{n} = \emptyset$, for every $I$, respectively. Furthermore, if $\gamma_{a_i}$ denotes the facet with perpendicular $a_i$, and $\gamma_I = \bigcap_{i \in I} \gamma_{a_i}$, then $f_I(x)= f_{\gamma_{I}}(x)$. 
\endproclaim

\proclaim{Lemma 4.4}
Let $f(x) \in \Cal{O}_K[x]$ be a   globally non-degenerate polynomial  with respect to its Newton
polyhedron $\Gamma(f)$,   $\gamma \subseteq \Gamma(f)$   a proper face,  and  $\Delta_{\gamma}$ its associated  cone. Let    $\{a_1,a_2,..,a_e\} \subset \Cal{D}$ be a  set of generators of  $\Delta_{\gamma}$. Then

$$\align 
Z_{\Delta_{\gamma}}( s,f, \chi )= 
\sum_{ y } A_y(q^{-s})Z(\Cal{O}_K ^{\times n }, s,f_y, \chi )+\\
\sum_{I \subseteq \{1,2,..,e\}}\frac{A_I(q ^{-s})Z(\Cal{O}_K ^{\times n }, s,f_I, \chi )}
{ \prod_{j \in I}(1- q ^{-\mid a_j\mid -m(a_j) s})},\\ 
\tag{4.13}
\endalign
$$
where $y$  runs through a finite number of points in $ \Bbb{N}^{n}$,  $A_y(q^{-s})$, $A_I(q^{-s}) \in \Bbb{Q}[q^{-s}]$,  with  $A_I(q^{-s})=\sum_k q^{-a_k(I)-b_k(I)s}$, $a_k(I), b_k(I) \in \Bbb{N}$, $f_y(x)$ and $f_I(x)$ are polynomials in $\Cal{O}_K[x]$ satisfying $Sing_{f_y}(K) \bigcap (K \smallsetminus \{0\})^{n} = \emptyset$, for every $y$, and   $Sing_{f_I}(K) \bigcap (K \smallsetminus \{0\})^{n} = \emptyset$, for every $I$, respectively. Furthermore, if $\gamma_{a_i}$ denotes the facet with perpendicular $a_i$, and $\gamma_I = \bigcap_{i \in I} \gamma_{a_i}$, then $\Gamma(f_I)= \gamma_{I}$.

\endproclaim
\demo{Proof}
By proposition 3.3 there exists a finite  partition of $\Delta_{\gamma}$ into  cones $\Delta_j$, such that  each  $\Delta_j$ is strictly positively spanned by  some    vectors from $\{a_1,a_2,..,a_e\}$ which are linearly independent over $\Bbb{R}$. Now, each cone $\Delta_j$ can be partitioned into a finite number of cones satisfying the conditions of corollary 4.3.
In order to verify this last assertion, we  observe that the set $\Delta_j \bigcap  \Bbb{N}^{n} $ admits the following partition:

$$\Delta_j \bigcap  \Bbb{N}^{n}=\left( \bigoplus_{i=1}^{e}a_i( \Bbb{N}\smallsetminus \{0\}\right) \bigcup \bigcup_{b} \left(b +   \bigoplus_{i=1}^{e}a_j \Bbb{N}\right), \tag{4.14}$$

where $b$ runs through a finite number of vectors in 

$$\Bbb{N}^{n} \bigcap \{\sum_{i=1}^{e}a_i\lambda_i \, \mid \, \lambda_i \in \Bbb{R}, \, 0 \leq \lambda_i<1, \, i=1,..,e\}.$$ 
Now the result follows from  corollary 4.3.\qed\quad
\enddemo

In the  proof of the above result, we did not use a partition of the cone $\Delta$ into  simple cones, because this approach produces a bigger list of candidates for the poles of $Z_{\Delta_\gamma}(s,f,\chi)$.

\subheading{ Proof of Theorem A}
(i) Given   a polynomial $f(x)\in \Cal{O}_K[x], f(0)= 0$, there exists  a partition of  $ \Bbb{R} _{+ } ^{n} $ of the form:

$$
\Bbb{R}_ + ^{n} = \{ (0,..,0)\}\bigcup \bigcup_{ \gamma}\Delta_{\gamma}, \tag{4.15}
$$
where $\gamma$ runs through all proper faces of $\Gamma(f)$, and  $\Delta_{\gamma}$ is a  cone strictly positive spanned by some vectors $a_1,..,a_e \in \Cal{D}$. In addition, $\Delta_{\gamma}$  is contained in an equivalence class of $\backsimeq$.
 From the above partition we obtain the following expression for $Z(s,f,\chi)$:
$$
Z(s,f,\chi) = \int_{\Cal{O}_K^{\times  n}}\chi(ac f(x))\mid f(x)\mid_K^{s}\mid dx\mid +
\sum_{\gamma}Z_{\Delta_{\gamma}}(s,f,\chi).
 \tag{4.16}
$$
 In (4.16) there are two different types of integrals: $Z(\Cal{O}_K^{\times  n},s,f,\chi)$,  and $Z_{\Delta_{\gamma}}(s,f,\chi)$.The integrals of the first type are rational functions of $q^{-s}$ with  poles satisfying  $ Re(s)=-1$ (cf. lemma 2.4). The second type of integrals are rational functions of $q^{-s}$ with  poles satisfying  condition (i) in the statement of theorem A (cf. lemma 4.4).

(ii) If $\chi \ne \chi_{\text{triv}}$, from (4.16) and lemma  2.4 follow  that $Z(s,f,\chi)$ is  equal to a polynomial, with degree bounded by a constant independent of $\chi$,  plus a finite sum  of functions  of the form 
$$
\frac{A_I(q ^{-s})Z(\Cal{O}_K ^{\times n }, s,f_I, \chi )}
{ \prod_{j \in I}(1- q ^{-\mid a_j\mid -m(a_j) s})},
\tag{4.17}
$$
where $f_I(x)$ denotes the restriction of $f(x)$ to the face $\gamma_I = \bigcap_{i \in I} \gamma_{a_i}$, and $\gamma_{a_i}$ denotes the facet with perpendicular $a_i$. The second part of the theorem follows from (4.17) by the following fact: 
$$
Z(\Cal{O}_K ^{\times n }, s,f_I, \chi )=0, \,\, \text{if the order of } \,\,\chi \,\, \text{does not divide some} \,\, m(a_j)\ne 0,\, j \in I.\tag{4.18}
$$
If the order of $\chi$ does not divide $m(a_j)$, with $a_j= (a_{1,j},a_{2,j},..,a_{n,j})$, then there exists an $u \in \Cal{O}_K^{\times}$ such that
$$ \chi^{m(a_j)}(u) \ne 1.\tag{4.19}$$
We set
$$
\CD
\phi_u: \Cal{O}_K^{\times n} @>>>  \Cal{O}_K^{\times n}\\
(x_1,x_2,..,x_n) @>>>(x_1 u^{a_{1,j}},x_2 u^{a_{2,j}},..,x_n u^{a_{n,j}}).\\
\endCD
\tag{4.20}
$$
The map   $\phi_u$ establishes a bijection  of $\Cal{O}_K^{\times n}$ to itself that preserves the Haar measure. 
By using (4.20) as change of variables in the integral $Z(\Cal{O}_K ^{\times n }, s,f_I, \chi )$,
it verifies that 
$$(1-\chi^{m(a_j)}(u) )Z(\Cal{O}_K ^{\times n }, s,f_I, \chi )=0.$$
Therefore, (4.19) implies $Z(\Cal{O}_K ^{\times n }, s,f_I, \chi )=0$. \qed\quad

\head
 {\bf 5. The largest pole of $Z(s,f,\chi_{\text{triv}})$}
\endhead
In this section we  prove theorem B. Its proof  will be accomplished by means of three preliminary results.

For a polynomial $f(x) \in \Cal{O}_K[x]$ globally non-degenerate  with respect to its Newton 
polyhedron $\Gamma(f)$, we set
$$
\beta(f) := \max_{\tau _j}\{-\frac{\mid a_j\mid}{m(a_j)}\},
$$
where $\tau _j$ runs through  all facets of  $\Gamma(f)$ satisfying $m(a_j) \ne 0$.  The point
 $$T_0 = (-\beta(f)^{-1},...,-\beta(f)^{-1})\in \Bbb{Q}^{n}$$ is the intersection point
 of the boundary of the Newton polyhedron $\Gamma(f)$ with the diagonal
 $\Delta = \{ (t,..,t) \; \mid \, t \in \Bbb{R} \}$ in $\Bbb{R}^{n}$. Let $\tau_0$ be 
the  face of smallest dimension of  $\Gamma(f)$  containing  $T_0$, and  $\rho$ its codimension, i.e. $\rho=$dim $\Delta_{\tau_0}$.

\proclaim{Proposition 5.1}
Let $f(x) \in \Cal{O}_K[x]$ be  a globally  non-degenerate polynomial  with respect to 
its Newton polyhedron $\Gamma(f)$. If $\beta(f) > -1$, then $\beta(f)$ is a pole of $Z(s,f,\chi_{\text{triv}})$  and its multiplicity    is   equal to  $\rho$.
\endproclaim

\demo{Proof}
First, we note that the multiplicity of the possible pole $\beta(f)$ is less then or equal to dim $\Delta_{\tau_0}= \rho$ (cf. formulas (4.16), (4.13), (2.7)). 
In order to prove that $\beta(f)$ is a pole of $Z(s,f,\chi_{\text{triv}}) $,  it is sufficient to show that
$$
\lim_{s \to \beta(f)}\left(1- q^{\beta(f)- s} \right)^{\rho}
Z(s,f,\chi_{\text{triv}}) > 0.
 \tag{5.1}
$$
 This last  assertion is a consequence of the  following  result (cf. (4.16), (4.13)):

\noindent {\bf Claim A} 

(i)
$$
 \lim_{s \to \beta(f)}\left(1- q^{\beta(f) - s} \right)^{\rho} 
\left(\frac{A_{I}(q ^{-s})Z(\Cal{O}_K ^{\times n }, s,f_I, \chi_{\text{triv}} )}
{ \Pi_{j \in I}(1- q ^{-\mid a_j\mid -m(a_j) s})}\right) \geq 0,
\tag{5.2}
$$ 

\noindent for every cone $\Delta_{\gamma}=\{ \sum_{i=1}^{e}a_i y_i \,\mid \, y_i \geq 0, \,\, \text{for all } \,
 i\}$, and every $I \subseteq \{1,2,..,e \}$.

(ii) There is a cone $\Delta_0$ and a subset 
$I_0$ of generators of  this cone such that  inequality  (5.2) is  strict.

The first part of the previous claim  follows from the following two facts. The first fact is
   
$$
 \lim_{s \to \beta(f)}
 \left(A_I(q ^{-s})Z(\Cal{O}_K ^{\times n }, s,f_I, \chi_{\text{triv}} ) \right)>0.
\tag{5.3}
$$
Since $A_I(q^{-s})=\sum_kq^{a_k(I)-b_k(I)s}$, with $a_k(I), b_k(I) \in \Bbb{N}$, inequality (5.3) follows from  noticing that
$$\lim_{s \to \beta(f)}\left( \frac{(1- q^{-1})q^{-s}}{1-q^{-1-s}}\right)>0, \,\, \text{when} \,\, \beta(f) >-1.$$

The  second fact is

$$\lim_{s \to \beta(f)}\left(1- q^{\beta(f)- s} \right)^{\rho} \left(\frac{1}
{ \prod_{j \in I}(1- q ^{-\mid a_j\mid -m(a_j) s})}\right) \geq 0. \tag{5.4}$$

The second part of the claim follows from the following reasoning. Let  $a_1, a_2,..,a_e$ be the unique primitive vectors perpendicular to the facets which contain $\tau_0$. There exists  a cone $\Delta_0$ in the partition into simplicial cones of $\Delta_{\tau_0}$ given by proposition 3.3. and $I_0 \subseteq \{1,2,..,e\}$ such that $\{a_i \mid i\in I_0\}$ is a set of $\rho$ linearly independent generators of $\Delta_0$, because the dimension of $\Delta_{\tau_0}$ is $\rho$. Then  inequality  (5.2) is strict  for the cone $\Delta_0$ and $I_0$. Thus, $\beta(f)$ is a pole of $Z(s,f,\chi_{\text{triv}})$ of multiplicity $\rho$.\qed\quad
\enddemo
\proclaim{Proposition 5.2}
Let $f(x) \in \Cal{O}_K[x]$ be a globally  non-degenerate  polynomial with respect to its Newton polyhedron $\Gamma(f)$, and $\gamma \subseteq \Gamma(f)$ a proper face. If $\sigma (\bar{f_{\gamma}},\Cal{O}^{\times \, n}_K)=\sigma (\bar{f_{\gamma}},\Cal{O}^{\times \, n}_K, \chi_{\text{triv}}) >0 $ then
$$
\lim_{s \to -1}\left(1- q^{-1- s} \right)
Z(\Cal{O}^{\times \, n}_K,s,f_{\gamma},\chi_{\text{triv}}) \ne 0. \tag{5.5}
$$
\endproclaim
\demo{Proof}
By using expansion (2.7), with $D=\Cal{O}^{\times \, n}_K$, and $m= C(f_{\gamma},\Cal{O}^{\times \, n}_K)+1$, we have that 
$$
\lim_{s \to -1}\left(1- q^{-1- s} \right)Z(\Cal{O}^{\times \, n}_K,s,f_{\gamma},\chi_{\text{triv}})=(q-1)\sigma (\bar{f_{\gamma}},\Cal{O}^{\times \, n}_K, \chi_{\text{triv}})+$$
$$
(q-1) \sum_{k = 1}^{m}q^{-kn}\left( 
\sum_{(P_1,...P_k) \in I_k}\sigma(\overline{f_{\gamma}}_{P_1,..,P_k},\chi_{\text{triv}})q^{E(P_1,..,P_k)}\right).
\tag{5.6}$$

Since the right side of  (5.6) is a sum of positive numbers, the result follows from the  hypothesis  
$ \sigma (\bar{f_{\gamma}},\Cal{O}^{\times \, n}_K, \chi_{\text{triv}})> 0 $.  \qed\quad
\enddemo

\proclaim{Proposition 5.3}
Let $f(x) \in \Cal{O}_K[x]$ be  a globally  non-degenerate polynomial  with respect to 
its Newton polyhedron $\Gamma(f)$. Let  $a_1, a_2,..,a_e$ be the unique primitive vectors perpendicular to the facets which contain $\tau_0$. If $\beta(f) =-1$, then $\beta(f)$ is a pole of $Z(s,f,\chi_{\text{triv}})$  with  multiplicity      less than or equal to $\rho +1$. Furthermore,  if  every face  $\gamma \supseteqq \tau_0$ satisfies $\sigma (\bar{f_{\gamma}},\Cal{O}^{\times \, n}_K)>0$,
then the multiplicity of the pole $\beta(f)$ is $\rho +1$.
\endproclaim

\demo{Proof}
 In the case $\beta(f)=-1$ the multiplicity of the possible pole $\beta(f)$ is less than or equal to $\rho +1$ because $Z(\Cal{O}_K ^{\times n }, s,f_I, \chi_{\text{triv}})$ may have a pole at $s=-1$ (cf. formulas (4.16), (4.13), (2.7)).
As in the case $\beta(f) > -1$, the result follows from  inequality (5.1) by claim A. In the case $\beta(f)=-1$, we may suppose that 
$$Z(\Cal{O}_K ^{\times n }, s,f_I, \chi_{\text{triv}} )= \frac{c_I (q^{-s})}{(1-q^{-1-s})},\tag{5.9}$$
\noindent where  $c_I(q^{-s})$ is a polynomial with positive coefficients (cf. expansion 2.7). The  proof of  claim A  involves the same ideas as in the  case $\beta(f) >-1$. 

The second  part of the proposition  is proved  as follows.   There exists a simplicial cone  $\Delta_0 \subseteq \Delta_{\tau_0}$  with dim $\Delta_0 =\rho$ (cf. final part of the proof of proposition 5.1). Let   $I_0$ be  a set of $\rho$ linearly independent generators of  $\Delta_0$. By   duality  this cone corresponds  to a face   $\gamma \supseteqq \tau_0$, and   $Z(\Cal{O}_K ^{\times n }, s,f_{I_0}, \chi_{\text{triv}})$  has a pole of multiplicity 1 at $s=-1$ (cf. proposition 5.2),  thus
$$
\lim_{s \to -1}\left(1- q^{-1 - s} \right)^{\rho +1} 
\left(\frac{A_{I_0}(q ^{-s})Z(\Cal{O}_K ^{\times n }, s,f_{I_0}, \chi_{\text{triv}} )}
{ \Pi_{j \in I_0}(1- q ^{-\mid a_j\mid -m(a_j) s})}\right) > 0.
\tag{5.10}
$$ \qed\quad
\enddemo
\subheading{ Proof of Theorem B}
The theorem follows from proposition 5.1 and proposition 5.3.
\qed\quad

\head
{\bf 6. Exponential sums}
\endhead

Let $\Psi$ be an additive character trivial  on  $\Cal{O}_K$ but 
not on $\Cal{P}_K^{-1}$. A such character is named {\it standard}. We  put  $z =u \pi ^{-m}$, $m \in \Bbb{N}\smallsetminus 
\{0\}$, $u \in \Cal{O}_K^{\times}$. To these data one associates the following exponential sum:
$$E(z,K,f) = q^{-nm} \sum_{x \, \text{mod} \, \Cal{P}_K^{m}}\Psi(uf(x)/\pi^{m}).$$
The  following corollary follows  theorem A,  theorem B above, and  proposition 1.4.5 of  [D-2],  
by writing $Z(s,f,\chi)$ in partial fractions.  
 
\proclaim{Corollary 6.1}
 (i) Let $f(x) \in \Cal{O}_K[x]$ be  a  globally  non-degenerate polynomial 
with respect to its Newton polyhedron $\Gamma(f)$, then
for $\mid z \mid$ big enough $E(z, K,f)$ is a finite $\Bbb{C}$-linear combination of functions of the form
$$\chi (ac (z)) \mid z \mid_K ^{\lambda} (log_q (|z|_K) )^{\beta},$$ with
 coefficients  independent of $z$, and with $\lambda \in \Bbb{C}$ a pole of
 $(1-q^{-1-s}) Z(s,f,\chi_{\text{triv}})$ or of $Z(s,f,\chi)$, 
$\chi \ne \chi_{\text{triv}}$, and $\beta \in \Bbb{N}$, $\beta \leqq $(multiplicity of pole 
$\lambda$) -1. Moreover all poles $\lambda$ appear effectively in this linear combination.

\noindent (ii) Let  $L$ be a global field,  and let $f(x) \in L[x]$ be  a 
globally non-degenerate  polynomial with respect to its Newton polyhedron $\Gamma(f)$, and suppose that $\beta(f) >-1$. 
For almost all non-archimedean  completions $L_v$ of $L$,  there exists a constant  $ C(L_v) \in \Bbb{R}$ satisfying
 
$$| E(z,L_v,f) | \leqq C(L_v) \mid z\mid_{L_v} ^{-\beta(f)}log_q (|z|_{L_v}) ^{\rho -1},  
 \,\,\ \text{for all}\,\, z \in L_v.$$

\endproclaim 
Igusa has conjectured that   $C(L_v)=1$ for almost all $v$ [I2]. 
This conjecture was proved by Denef and   Sperber when $K$ has characteristic zero,  $f$ is 
a non-degenerate polynomial, and  the face of the
 Newton polyhedron  which cuts the diagonal  does not have  vertex
  in $\{0,1  \}^{n}$ [D-Sp]. Corollary 6.1  permits us to extent the result of Denef and Sperber to positive characteristic using the methods in [D-Sp].

\head
{\bf 7. Examples}
\endhead
\example{Example 7.1}
In this example,  we compute $Z(s,f,\chi_{triv})=Z(s,f)$, for $f(x,y) =
x^{2}+xy+y^{2}$, when the characteristic  of $K$ is different  from 2, 3, and analyze the behavior of the pole $s=-1$.  In this case $Sing_{f}(K) =\{(0,0)\}$, and the Newton polygon has only a compact segment with supporting hyperplane  $ x + y=2$. The polynomial $f$ is globally non-degenerate with respect to its Newton polygon. 

One easily verifies  that $\Bbb{R}_+^{2} \smallsetminus \{(0,0)\}$ can be partitioned into equivalence classes modulo $ \backsimeq$, as follows.

If
$$\Delta_1:=\{(0, a)\, \mid \,a > 0\},$$
$$\Delta_2:=\{(b, a+b)\, \mid \,a,b >0\},$$
$$\Delta_3:=\{(a, a)\, \mid \,a > 0\},$$
$$\Delta_4:=\{(a+b, a)\, \mid \,a,b >0\},$$
$$\Delta_5:=\{(a, 0)\, \mid \,a > 0\},$$

then
$$\Bbb{R}_+^{2}= \{(0,0)\} \bigcup \bigcup_{i=1}^{5}\Delta_i,$$
\noindent and
$$Z(s,f)= Z(\Cal{O}_K^{\times 2}, s,f)+\sum_{i=1}^{5}Z_{\Delta_i}(s,f).$$

{\bf Calculation of $Z(\Cal{O}_K^{\times 2}, s,f)$, and $Z_{\Delta_1}(s,f)$}

By using the stationary phase formula, we obtain
$$
Z(\Cal{O}_K^{\times 2}, s,f)= \nu(\overline{f},\Cal{O}_K^{\times 2}) +\sigma(\overline{f},\Cal{O}_K^{\times 2})
\frac{(1-q^{-1})q^{-1}}{(1-q^{-1-s})}.
\tag{7.1}
$$

On the other hand, it is simple to verify that    $Z_{\Delta_1}(s,f)=q^{-1}(1-q^{-1})$.

{\bf Calculation of $Z_{\Delta_2}(s,f)$ and $Z_{\Delta_3}(s,f)$}

$$
Z_{\Delta_2}(s,f)= \sum_{a, b = 1}^{\infty}q^{-a-2b}\int_{\Cal{O}_K^{\times 2}}\mid  
  \pi^{2b}x^{2}+ \pi^{a+2b}xy +\pi^{2a+2b}y^{2}  \mid_K^{s} \mid dx dy\mid
$$
$$=\frac{q^{-3-2s}(1-q^{-1})}{(1-q^{-1-s})(1+q^{-1-s})}.
\tag{7.2}
$$

$$\align
Z_{\Delta_3}(s,f)=\sum_{a\geq 1}^{\infty}q^{-2a} \int_{\Cal{O}_K^{\times 2}}\mid \pi^{2a}x^{2} + \pi^{2a}xy + \pi^{2a}y^{2}\mid_K^{s} \mid dx dy \mid=\\
\frac{q^{-2-2s}}{(1- q^{-1-s})(1+ q^{-1-s})}\left(\nu(\overline{f},\Cal{O}_K^{\times 2}) + \sigma(\overline{f},\Cal{O}_K^{\times 2}) \frac{(1- q^{-1})q^{-s}}{(1-q^{-1-s})}\right).\\
\tag{7.3}
\endalign
$$

{\bf Calculation of $Z_{\Delta_4}(s,f)$ and $Z_{\Delta_5}(s,f)$}

$$
Z_{\Delta_4}(s,f)= \sum_{a,b \geq 1}^{\infty}q^{-2a-b} \int_{\Cal{O}_K^{\times 2}}\mid \pi^{2a+2b}x^{2} + \pi^{2a+b}xy + \pi^{2a}y^{2}\mid_K^{s} \mid dx dy \mid=
$$
$$
\frac{q^{-3-2s}(1- q^{-1})}{(1- q^{-1-s})(1+ q^{-1-s})}.
\tag{7.4}
$$

$$
Z_{\Delta_5}(s,f)=q^{-1}(1-q^{-1}). \tag{7.5}
$$

From the above calculations, we obtain
$$
\lim_{s \to -1}\left(1- q^{-1- s} \right)^{2}
Z(s,f)= \frac{\sigma(\overline{f},\Cal{O}_K^{\times 2}) (q-1)}{2}.\tag{7.6}$$
Now suppose that  $K=\Bbb{Q}_p$,  with $p\ne 2,3$. Since
$$\sigma(f,\Cal{O}^{\times \, 2}_K) =p^{2}\text{Card}\left( \{ (u,v)\in \Bbb{F}^{\times \, 2}_p \mid 
\overline{f}(u,v)=0\}\right)=\cases
0, \,\, \text{if} \,\, p\equiv 5,11 \,\ \text{mod} \,\ 12,\\
2p^{-2}(p-1), \,\, \text{if} \,\, p\equiv 1,7 \,\ \text{mod} \,\ 12,\\
\endcases
$$

it follows from  (7.6) that
$$
\lim_{s \to -1}\left(1- p^{-1- s} \right)^{2}
Z(s,f)= \cases
0, \,\, \text{if} \,\, p\equiv 5,11 \,\ \text{mod} \,\ 12,\\
p^{-2}(p-1)^{2}, \,\, \text{if} \,\, p\equiv 1,7 \,\ \text{mod} \,\ 12.\\
\endcases
\tag{7.7}
$$

Thus  $Z(s,f)$ has a pole  at $s=-1$ of multiplicity $\rho +1 =2$, when  $$\text{Card}\left( \{ (u,v)\in \Bbb{F}^{\times \, 2}_p \mid \overline{f_{\tau_0}}(u,v)=0\}\right)=\text{Card}\left( \{ (u,v)\in \Bbb{F}^{\times \, 2}_p \mid \overline{f}(u,v)=0\}\right)>0.$$ Otherwise the multiplicity  is $\rho=1$.   
\endexample

\example {Example 7.2}
In this example, by using the method of lemma 4.4,  we compute the local zeta function attached to the polynomial $f(x,y)= x^{2}y^{2} + x^{5} +y^{5} \in K[x,y]$, when the characteristic of $K$ is different from 2, 5. This polynomial is globally non-degenerate with respect to its Newton polyhedron. 

One easily verifies  that $\Bbb{R}_+^{2} \smallsetminus \{(0,0)\}$ can be partitioned into equivalence classes modulo $ \backsimeq$, as follows.

If

$$\Delta_1:=\{(0, a)\, \mid \,a > 0\},$$
$$\Delta_2:=\{(2b, a+3b)\, \mid \,a,b >0\},$$
$$\Delta_3:=\{(2a, 3a)\, \mid \,a > 0\},$$
$$\Delta_4:=\{(2a+3b, 3a+2b)\, \mid \,a,b >0\},$$
$$\Delta_5:=\{(3a, 2a)\, \mid \,a > 0\},$$
$$\Delta_6:=\{(3a+b, 2a)\, \mid \,a,b > 0\},$$
$$\Delta_7:=\{(a, 0)\, \mid \,a > 0\},$$
then 
$$\Bbb{R}_+^{2}= \{(0,0)\} \bigcup \bigcup_{i=1}^{7}\Delta_i,$$
where each $\Delta_i$ is exactly an equivalence class modulo $ \backsimeq$.

{\bf Calculation of $Z(\Cal{O}_K^{\times 2}, s,f)$, and $Z_{\Delta_1}(s,f)$}

By using the stationary phase formula, we obtain
$$
Z(\Cal{O}_K^{\times 2}, s,f)= \nu(\overline{f},\Cal{O}_K^{\times 2}) +\sigma(\overline{f},\Cal{O}_K^{\times 2})\frac{(1-q^{-1-s})}{1-q^{-1-s}}.
\tag{7.8}
$$
On the other hand, it is simple to verify that    $Z_{\Delta_1}(s,f)=q^{-1}(1-q^{-1})$.

{\bf Calculation of $Z_{\Delta_2}(s,f)$ and $Z_{\Delta_3}(s,f)$}

The cone $\Delta_2$ is not a simple. In this case, one verifies  that there is only one element  in $\Delta_2 \bigcap \Bbb{N}^{2}$ satisfying $0 \leq a <1$, $0\leq b <1$. This element is $(1,2)= (0,1)\frac{1}{2} +(2,3)\frac{1}{2}$. Thus
$$\Delta_2 \bigcap \Bbb{N}^{2} = \{ (0,1) (\Bbb{N} \smallsetminus \{0\})+ (2,3)(\Bbb{N} \smallsetminus \{0\})\} \bigcup \{(1,2)+ (0,1) \Bbb{N}+ (2,3)\Bbb{N}\}. \tag{7.9}
$$
From the partition (7.9), we obtain that
$$\align
Z_{\Delta_2}(s,f)= \sum_{a, b = 1}^{\infty}q^{-a-5b}\int_{\Cal{O}_K^{\times 2}}\mid  
  \pi^{2a+10b}x^{2}y^{2} +\pi^{10b}x^{5} + \pi^{5a+15b}y^{5}  \mid_K^{s} \mid dx dy\mid\\
+\sum_{a, b = 0}^{\infty}q^{-a-5b-3 }\int_{\Cal{O}_K^{\times 2}}\mid  
  \pi^{2a+10b+6}x^{2}y^{2} +\pi^{10b+5}x^{5} + \pi^{5a+15b+10}y^{5}  \mid_K^{s} \mid dx dy\mid=\\
\frac{q^{-5-10s}}{1-q^{-5-10s}}q^{-1}(1-q^{-1})+\frac{q^{-3-5s}}{1-q^{-5-10s}}(1-q^{-1})
=\frac{(1-q^{-1})(q^{-3-5s}+ q^{-6-10s})}{1-q^{-5-10s}}.\\
\tag{7.10}
\endalign
$$
By applying proposition 4.1, and then the stationary phase formula to  $Z_{\Delta_3}(s,f)$, one obtains
$$\align
Z_{\Delta_3}(s,f)=\sum_{a = 1}^{\infty}q^{-5a-10as} \int_{\Cal{O}_K^{\times 2}}\mid y^{2} + x^{3}\mid_K^{s} \mid dx dy \mid=\\
\frac{q^{-5-10s}}{1- q^{-5-10s}}\left(\nu(\overline{f},\Cal{O}_K^{\times 2}) + \sigma( \overline{f},\Cal{O}_K^{\times 2}) \frac{(1- q^{-1})q^{-s}}{(1-q^{-1-s})}\right).\\
\tag{7.11}
\endalign
$$

{\bf Calculation of $Z_{\Delta_4}(s,f)$ and $Z_{\Delta_5}(s,f)$}

The cone $\Delta_4$ is not a simple, thus we proceed as in the computation of $Z_{\Delta_2}(s,f)$, i.e. we  find $ 0 \leq a <1$, $0 \leq  b <1$, such that 
$$ (2,3)a + (3,2)b \in \Bbb{N}^{2} \bigcap \Delta_4.$$

If $a=b$, one  finds   immediately that  $(2,3)\frac{i}{5}+(3,2)\frac{i}{5} \in \Bbb{N}^{2} \bigcap \Delta_4$, $i=1,2,3,4$. The case $  a\ne b$ cannot occur. Suppose that $(m,n)\in 
\Bbb{N}^{2} \bigcap \Delta_4$, with  $b >a$, $a\ne 0$, $b \ne 0,$ ($a=0$ or $b=0$ cannot occur), i.e.
$$
m= 2a+3b, \,\,\, n=3a+2b, \,\, m,n \in \Bbb{N} \smallsetminus \{0\}, \,\, 0< a < b < 1.\tag{7.12}$$
From (7.12), we get $b-a =m-n$,  but this is impossible because $0 < b-a < 1$, and $m-n \geq 1$. If $a >b$ then $a-b=n-m$ and the same argument applies.

Therefore, we have the following partition for $\Bbb{N}^{2} \bigcap \Delta_4$:
$$\align
\Bbb{N}^{2} \bigcap \Delta_4= \{(2,3) (\Bbb{N}\smallsetminus \{0\}) + (3,2) (\Bbb{N}\smallsetminus \{0\})\} \bigcup \bigcup_{i=1}^{4} \{(i,i) +(2,3) \Bbb{N} + (3,2) \Bbb{N}\}.\tag{7.13}
\endalign
$$

From the partition (7.13), we obtain that
$$
Z_{\Delta_4}(s,f)= \left( \frac{(1-q^{-1})(q^{-5-10s})}{1-q^{-5-10s}} \right) ^{2} +
\left( \sum_{i=1}^{4} q ^{-2i-4is}\right) 
\left( \frac{1-q^{-1}}{1-q^{-5-10s}} \right) ^{2}. 
\tag{7.14}
$$
For $Z_{\Delta_5}(s,f)$, we get

$$
Z_{\Delta_5}(s,f)= \frac{q^{-5-10s}}{1-q^{-5-10s}}\left(
 \nu(\overline{f},\Cal{O}_K^{\times 2}) + \sigma(\overline{f},\Cal{O}_K^{\times 2}) \frac{(1- q^{-1})q^{-s}}{1-q^{-1-s}} \right).\tag{7.15}
$$

{\bf Calculation of $Z_{\Delta_6}(s,f)$}

In the computation of the integral $Z_{\Delta_6}(s,f)$, we use the following partition:
$$\Delta_6 \bigcap \Bbb{N}^{2} = \{ (3,2) (\Bbb{N}\smallsetminus \{0\})+ (1,0)(\Bbb{N}\smallsetminus \{0\})\} \bigcup  \{(2,1) + (3,2)\Bbb{N}+ (1,0)\Bbb{N}\}.\tag{7.16}$$
From the above partition, we get
$$Z_{\Delta_6}(s,f) =(1- q^{-1}) \frac{q^{-3-5s}+ q^{-6-10s}}{1-q^{-5-10s}}.\tag{7.17}$$

{\bf Calculation of $Z_{\Delta_7}(s,f)$}
$$Z_{\Delta_7}(s,f)=q^{-1}(1-q^{-1}).\tag{7.18}$$

Now, with $\beta(f)=-\frac{1}{2}$, and $\rho =2$, it holds that  
$$
\lim_{s \to \beta(f)}\left(1- q^{\beta(f)- s} \right)^{\rho}
Z(s,f)=\lim_{s \to \beta(f)}\left(1- q^{\beta(f)- s} \right)^{\rho}
Z_{\Delta_4}(s,f)=\frac{(1-q^{-1})^{2}}{50}.
$$ 

\endexample

\refstyle{C}

\enddocument